\documentclass{article}
\usepackage{amsmath}
\usepackage{amsfonts}

\newenvironment{proof}[1][Proof]{\noindent\textbf{#1.} }{\ \rule{0.5em}{0.5em}}
\input{tcilatex}

\begin{document}

\title{General proof of a limit related to AR(k) model of Statistics}
\date{}
\author{Jan \ VRBIK \\
Department of Mathematics\\
Brock University, 500 GLenridge Ave.\\
St. Catharines, Ontario, Canada, L2S 3A1}
\maketitle

\begin{abstract}
Computing moments of various parameter estimators related to an
autoregressive model of Statistics, one needs to evaluate several
non-trivial limits. This was done by \cite{lui} for the case of two, three
and four dimensions; in this article, we present a proof of a fully general
formula, based on an ingenious solution of \cite{fedor}.
\end{abstract}

\section{Introduction}

The autoregressive model of Statistics generates a random sequence of
observations by%
\begin{equation}
X_{i}=\alpha _{1}X_{i-1}+\alpha _{2}X_{i-2}+...+\alpha
_{k}X_{i-k}+\varepsilon _{i}  \label{ARk}
\end{equation}%
where $\varepsilon _{i}$ are independent, Normally distributed random
variables with the mean of $0$ and the same standard deviation, and $k$ is a
fixed integer, usually quite small (e.g. $k=1$ defines the so called Markov
model). The sufficient and necessary condition for the resulting sequence to
be asymptotically stationary is that all $k$ solutions of the characteristic
polynomial%
\begin{equation}
\lambda ^{k}=\alpha _{1}\lambda ^{k-1}+\alpha _{2}\lambda ^{k-2}+...+\alpha
_{k}  \label{CHP}
\end{equation}%
are, in absolute value, smaller than $1$ (this is then assumed from now on).

The $j^{th}$-order serial correlation coefficient $\rho _{j}$ (between $%
X_{i} $ and $X_{i+j}$) is then computed by%
\begin{equation}
\rho _{j}=A_{1}\lambda _{1}^{|j|}+A_{2}\lambda _{2}^{|j|}+...+A_{k}\lambda
_{k}^{|j|}  \label{rho}
\end{equation}%
where the $\lambda _{i}$'s are the $k$ roots of (\ref{CHP}), and the $A_{i}$
coefficients are themselves simple functions of these roots. Note that the
absolute value of each root must be smaller than $1$ if the resulting
stochastic process is be stationary.

Computing the first few moments of various estimators (of the $\alpha _{i}$
parameters) boils down to computing moments of expressions of the 
\begin{equation}
\sum_{i=1}^{n}X_{i}  \label{SX}
\end{equation}%
and 
\begin{equation}
\sum_{i=1}^{n-j}X_{i}X_{i+j}  \label{DX}
\end{equation}%
type, where $X_{1},$ $X_{2},...X_{n}$ is a collection of $n$ consecutive
observations (assuming that the process has already reached its stationary
phase).

This in turn requires evaluating various summations (see \cite{vrbik}), of
which the most difficult has the form of%
\begin{equation}
\sum_{i_{1},i_{2},...i_{k}=1}^{\tilde{n}}\lambda
_{1}^{|i_{1}-i_{2}+s_{1}|}\lambda _{2}^{|i_{2}-i_{3}+s_{2}|}...\lambda
_{k}^{|i_{k}-i_{1}+s_{k}|}  \label{S4}
\end{equation}%
where $\lambda _{1},$ $\lambda _{2},...\lambda _{k}$ are the $\lambda _{i}$
roots (some may be multiple), $s_{1},$ $s_{2},...s_{k}$ are (small)
integers, and $\tilde{n}$ indicates that the upper limit equals to $n$,
adjusted in the manner of (\ref{DX}).

For small $k$, it is possible (but rather messy - see \cite{LIU}) to \emph{%
exactly} evaluate (\ref{S4}) and realize that the answer will always consist
of three parts:

\begin{itemize}
\item terms proportional to $\lambda _{i}^{n},$ which all tend to zero (as $%
n $ increases) `exponentially',

\item terms which stay constant as $n$ increases,

\item terms proportional to $n.$
\end{itemize}

Luckily, to build an approximation which is usually deemed sufficient (see 
\cite{vrbik}), we need to find only the $n$ proportional terms. These can be
extracted by dividing (\ref{S4}) by $n$ and taking the $n\rightarrow \infty $
limit. Incidentally, this results in the following (and most welcomed)
simplification: the corresponding answer will be the \emph{same} regardless
of the $\tilde{n}$ adjustments (thus, we may as well use $n$ instead), and
will similarly \emph{not} depend on the individual $s_{i}$'s, but only on
the \emph{absolute} value of their \emph{sum}, as the following statement
indicates.

\section{The main theorem}

\begin{eqnarray}
&&A\overset{\text{def}}{=}\lim_{n\rightarrow \infty }\frac{1}{n}%
\sum_{i_{1},i_{2},...i_{k}=1}^{n}\lambda _{1}^{|i_{1}-i_{2}-s_{1}|}\lambda
_{2}^{|i_{2}-i_{3}-s_{2}|}...\lambda _{k}^{|i_{k}-i_{1}-s_{k}|}  \label{A} \\
&=&\sum_{j=1}^{k}\lambda _{j}^{S+k-1}\prod\limits_{\substack{ \ell =1  \\ %
\ell \neq j}}^{k}\frac{1-\lambda _{\ell }^{2}}{(\lambda _{j}-\lambda _{\ell
})(1-\lambda _{j}\lambda _{\ell })}  \notag
\end{eqnarray}%
where $S=|s_{1}+s_{2}+...+s_{k}|.$

\begin{proof}
Define%
\begin{equation}
B_{S}\overset{\text{def}}{=}\sum_{m_{1}+m_{2}+...+m_{k}=S}\lambda
_{1}^{|m_{1}|}\lambda _{2}^{|m_{2}|}...\lambda _{k}^{|m_{k}|}  \label{B}
\end{equation}%
where $S$ is the non-negative integer of the theorem.

When $S\geq 0$, a term of $B_{S}$ and a term of the $A$ summation are
considered identical (we also say that they \emph{match} each other) only
when $m_{1}=-i_{1}+i_{2}+s_{1},$ $m_{2}=-i_{2}+i_{3}+s_{2},$ ... $%
m_{k-1}=-i_{k-1}+i_{k}+s_{k-1}$ (implying $m_{k}=-i_{k}+i_{1}+s_{k},$ since
the $m$'s and $s$'s must add up to the same $S,$ and the $\acute{\imath}$'s
cancel); note that this also implies (but not the reverse) that such
matching terms have the same \emph{value}. On the other hand, when $S<0$, we
declare them identical when $m_{1}=i_{1}-i_{2}-s_{1},$ etc. instead. From
now on, we assume that $S\geq 0$ to avoid a trivial duplication of all
subsequent arguments.

Clearly, each term of the $A$ summation matches a term of $B_{S}$: just take 
$m_{p}=-i_{p}+i_{p+1}+s_{p}$ where $p=1,$ $2,...k$, with the understanding
that $i_{k+1}=i_{1}$.

At the same time, no term of $B_{S}$ is matched by more than $n$ terms of
the $A$ summation, since once you select $i_{1}$ (from any of its $n$
possible values), all the remaining $i$'s are uniquely determined by $%
i_{2}=m_{1}+i_{1}-s_{1,}$ $i_{3}=m_{2}+i_{2}-s_{2,}$ etc., resulting in a
term of $A$ \emph{only} when all of these turn out to be between $1$ and $n$
(inclusive).

This proves that $A\leq B_{S}.$

Since $|\lambda _{1}|^{|m_{1}|}\leq 1,$ 
\begin{equation*}
B_{S}<\sum_{m_{2},m_{3},...m_{k}=-\infty }^{\infty }|\lambda
_{2}|^{|m_{2}|}|\lambda _{3}|^{|m_{3}|}...|\lambda _{k}|^{|m_{k}|}
\end{equation*}%
implying that the $B_{S}$ sum is (absolutely) convergent; let $B_{\infty }$
denote its actual value. This means that \emph{any} number smaller that $%
B_{\infty }$ (say $B_{0}$) can be \emph{exceeded} by a sum of \emph{finitely}
many terms of $B_{S}$ (this is true for any convergent series).

Now, let us go back to counting how many terms of the $A$ summation match a
single, specific term of $B_{S}$; we have already seen that, starting with
any one of the possible $n$ values of $i_{1},$ the subsequent $i$'s would be
computed by%
\begin{equation*}
i_{p}=i_{1}+\sum_{j=1}^{p-1}(m_{j}-s_{j})\ \ \ \ \ \ \ \ \ \text{where \ \ \ 
}p=2...k
\end{equation*}%
matching a term of the $A$ summation only when they are all in the $1$ to $n$
range, i.e. when 
\begin{equation*}
1\leq i_{1}+\min_{p=2...k}\sum_{j=1}^{p-1}(m_{j}-s_{j})
\end{equation*}%
and%
\begin{equation*}
i_{1}+\max_{p=2...k}\sum_{j=1}^{p-1}(m_{j}-s_{j})\leq n
\end{equation*}%
This implies that, for each choice of $i_{1}$ which meet%
\begin{equation*}
1-\min_{p=2..k}\sum_{j=1}^{p-1}(m_{j}-s_{j})\leq i_{1}\leq
n-\max_{p=2..k}\sum_{j=1}^{p-1}(m_{j}-s_{j})
\end{equation*}%
we get a legitimate term of the $A$ summation (matching and having the same
value as the specific term of $B_{S}$); we thus have%
\begin{equation*}
n-\max_{p=2...k}\sum_{j=1}^{p-1}(m_{j}-s_{j})+\min_{p=2..k}%
\sum_{j=1}^{p-1}(m_{j}-s_{j})
\end{equation*}%
such terms in total. Dividing their sum by $n$ and taking the $n\rightarrow
\infty $ limit thus yields the value of the specific $B_{S}$ term.

This can be repeated for any term of the finite sum of the previous
paragraph; thus we get $A\geq B_{0}.$ And, since we can make $B_{0}$ as
close to $B_{\infty }$ as we wish, this implies that $A\geq B_{\infty }.$

We have thus shown that (\ref{A}) and (\ref{B}) have the same value.

We now define the following Laurent series of the $B_{S}$ sequence (allowing 
$S$ to have \emph{any} integer value, and assuming that $\max_{\ell
=1...k}|\lambda _{\ell }|<|t|<\min_{\ell =1...k}|\lambda _{\ell }|^{-1}$),
namely%
\begin{eqnarray*}
&&F(t)\overset{\text{def}}{=}\sum_{S=-\infty }^{\infty
}t^{S}\sum_{m_{1}+m_{2}+...+m_{k}=S}\lambda _{1}^{|m_{1}|}\lambda
_{2}^{|m_{2}|}...\lambda _{k}^{|m_{k}|}\overset{}{=} \\
&&\sum_{m_{1},m_{2},...,m_{k}=-\infty }^{\infty }t^{m_{1}}\lambda
_{1}^{|m_{1}|}t^{m_{2}}\lambda _{2}^{|m_{2}|}...t^{m_{k}}\lambda
_{k}^{|m_{k}|}\overset{}{=} \\
&&\sum_{m_{1},m_{2},...,m_{k}=-\infty }^{\infty }\prod\limits_{\ell
=1}^{k}t^{m_{\ell }}\lambda _{\ell }{}^{|m_{\ell }|}\overset{}{=}%
\prod\limits_{\ell =1}^{k}\sum_{m=-\infty }^{\infty }t^{m}\lambda _{\ell
}{}^{|m|}\overset{}{=} \\
&&\prod\limits_{\ell =1}^{k}\left( \sum_{m=0}^{\infty }(t\lambda _{\ell
})^{m}+\sum_{m=-\infty }^{-1}t^{m}\lambda _{\ell }{}^{-m}\right) \overset{}{=%
}\prod\limits_{\ell =1}^{k}\left( \sum_{m=0}^{\infty }(t\lambda _{\ell
})^{m}+\sum_{m=1}^{\infty }(\frac{\lambda _{\ell }}{t})^{m}\right) \overset{}%
{=} \\
&&\overset{}{=}\prod\limits_{\ell =1}^{k}\left( \frac{1}{1-t~\lambda _{\ell }%
}+\frac{\frac{\lambda _{\ell }}{t}}{1-\frac{\lambda _{\ell }}{t}}\right) 
\overset{}{=}\sum\limits_{j=1}^{k}\left( \frac{C_{j}}{1-t~\lambda _{j}}+%
\frac{D_{j}}{t-\lambda _{j}}\right)
\end{eqnarray*}%
where the last expression is the partial-fraction expansion of the previous 
\emph{rational} function of $t$ (the roots of the common denominator are the 
$\lambda $'s and their inverses). We can now get a formula for $B_{\infty }$
(and thus for our $A$ limit) as a coefficient of $t^{S}$ of the last
expression. Since only the $C_{j}$ part contributes to non-negative powers
of $t,$ and%
\begin{eqnarray*}
&&C_{j}\overset{}{=}\left. F(t)(1-t~\lambda _{j})\right\vert _{t=\lambda
_{j}^{-1}}\overset{}{=}\prod\limits_{\substack{ \ell =1  \\ \ell \neq j}}%
^{k}\left( \frac{1}{1-\frac{\lambda _{\ell }}{\lambda _{j}}}+\frac{\lambda
_{j}\lambda _{\ell }}{1-\lambda _{j}\lambda _{\ell }}\right) \overset{}{=} \\
&&\overset{}{=}\lambda _{j}^{k-1}\prod\limits_{\substack{ \ell =1  \\ \ell
\neq j}}^{k}\left( \frac{1}{\lambda _{j}-\lambda _{\ell }}+\frac{\lambda
_{\ell }}{1-\lambda _{j}\lambda _{\ell }}\right) \overset{}{=}\lambda
_{j}^{k-1}\prod\limits_{\substack{ \ell =1  \\ \ell \neq j}}^{k}\frac{%
1-\lambda _{\ell }^{2}}{(\lambda _{j}-\lambda _{\ell })(1-\lambda
_{j}\lambda _{\ell })}
\end{eqnarray*}%
the final formula is therefore given by%
\begin{equation*}
\sum\limits_{j=1}^{k}\lambda _{j}^{S}C_{j}
\end{equation*}%
(note that the coefficient of $t^{S}$ in the expansion of $(1-t~\lambda
_{j})^{-1}$ is $\lambda _{j}^{S}$).

This proves the original statement.
\end{proof}

\section{Conclusion}

The formula of (\ref{A}) then enables us to evaluate all the expected values
needed to deal with any autoregressive model of type (\ref{ARk}). Note that
in some cases the set of $\lambda _{i}$ values may consist of only a \emph{%
subset} of of roots of (\ref{CHP}); this only reduces the value of $k$ and
makes the result that much easier.

A modification of the formula is needed when some of the $\lambda _{i}$'s
are \emph{identical}; in that case all we have to do is to evaluate the
formula's corresponding limit, such as $\lambda _{5}\rightarrow \lambda _{2}$
when the two $\lambda $'s have the same value (in the case of triple roots,
we would need to take two consecutive limits, etc.). This yields a multitude
of new (and rather messy) formulas not worth quoting - suffices to say that
they all result (as they must) in a finite expression.

A further challenge would be to find the \emph{constant} part of (\ref{S4}).

\end{document}